%% file: cycling_project.tex
\documentclass[twocolumn]{svjour_blank}

\usepackage{mathptmx} 

\usepackage[utf8]{inputenc}

\usepackage{mathtools, cuted}
 \usepackage{amsmath}
 \usepackage{amsfonts}
 \usepackage{amssymb}
\usepackage{graphicx}
\graphicspath{{include/}}
\usepackage{subfig}
\usepackage[hyphens]{url} 
\usepackage{color}
\usepackage{mathtools}
\usepackage{subfiles}
\usepackage{todonotes}
\usepackage{overpic}
\usepackage[normalem]{ulem}
\usepackage{comment}
\definecolor{deeppurple}{RGB}{102, 0, 153}

\newcommand{\maxm}{\mathrm{max}}
\newcommand{\minm}{\mathrm{min}}
\newcommand{\re}{\mathrm{e}}
\newcommand{\lr}[1]{\left(#1\right)}

\newcommand{\rd}{\mathrm{d}}
\newcommand{\tf}{t_f}
\newcommand{\ta}{t_a}
\newcommand{\tp}{t_p}
\newcommand{\Deltat}{\Delta t}

\newcommand{\ian}[1]{\todo[inline,color=red!40]{Ian: #1}}

\newcommand{\jcv}[1]{{\color{teal} #1}}
\newcommand{\jcvc}[1]{{\color{black} #1}}

\begin{document}

\title{A mathematical model for optimal breakaways in cycling: balancing energy expenditure and crash risk.}


\author{J.~Chico-V\'{a}zquez$^1$ \and I.~M.~Griffiths$^{1\,\ast}$}

\institute{$^1$Mathematical Institute, University of Oxford, Oxford, OX2 6GG \and $^\ast$ian.griffiths@maths.ox.ac.uk}

\maketitle

\vspace{-50mm}
\begin{abstract}
We present a mathematical model for optimizing breakaway strategies in competitive cycling, balancing power expenditure, aerodynamic drag, and crashing. Our  framework incorporates probabilistic crash dynamics, allowing a cyclist's risk tolerance to shape optimal tactics. We define an objective function that accounts for both  finish time differences and the probability of crashing, which we optimize subject to an energy expenditure constraint. We demonstrate the methodology for a flat stage with a simple constant-power breakaway. We then extend this analysis to account for fatigue-driven power decay, and varying terrain and race conditions. We highlight the importance of strategy by demonstrating that carefully planned decision making can lead to a race win even when the energy expenditure is low. Our results highlight and quantify the fact that, at the elite level, success often depends as much on minimizing risk as on maximizing physical output. 


\color{black}


\end{abstract}


\section{Introduction} \label{sec:intro}

In competitive cycling, cyclists typical ride for large portions of the race as a group, or \emph{peloton}. This has the advantage of significant reductions in aerodynamic drag. However, to win a race, a cyclist must choose to break away from the peloton at some strategic position during the race. To execute a successful breakaway, cyclists must balance energy expenditure, aerodynamic drag, and the inherent risk of crashing, often under complex and dynamic race conditions. To date, existing research has primarily focused on isolated aspects of racing, such as drag reduction through drafting~\cite{blockenCFDSimulationsAerodynamic2013}, nutrition or the biomechanics of power generation~\cite{atkinsonScienceCyclingCurrent2003,martinValidationMathematicalModel1998}. However, exploration of the specific interplay between these factors in breakaway situations is limited. Understanding this interaction is essential for optimizing race strategy, particularly for athletes seeking to achieve a tactical advantage.

Aerodynamic drag plays a significant role in determining the energy expenditure of a cyclist, especially when drafting within the peloton or executing a breakaway. As previous studies have shown, aerodynamic resistance accounts for up to 90\% of the overall resistance a cyclist experiences at race speeds~\cite{vandruenenAerodynamicAnalysisUphill2021,kyleReductionWindResistance1979}. While drafting can significantly reduce drag -- by as much as 80\% for cyclists in a peloton \cite{oldsMathematicsBreakingAway1998,blockenAerodynamicDragCycling2018} -- a breakaway rider is faced with the challenge of overcoming this loss of aerodynamic advantage. 

Gaul \emph{et al.}~\cite{gaulOptimizingBreakawayPosition2018} considered the problem of breakaway optimization by focusing on cyclist exhaustion and stamina, but did not consider the peloton structure or the important role crashing plays in their optimization strategy.

Crashing also contributes substantially to strategic \linebreak decision-making. Even a minor crash with no injuries can mean a rider loses around a minute, in addition to the extra energy expenditure required to catch up with the peloton. In a sport where victories are achieved with razor-thin margins, a single crash can make or break a contender's attempt. Moreover, a single major crash can lead to severe injuries, and ultimately a forcing a contender to abandon  an event. 

Crashes are influenced by environmental variables, race conditions and cyclist behaviour~\cite{dashFactorsImpactingBike2022,salmonBicycleCrashContributory2022}. In peloton settings, a higher density of riders increases the likelihood of collisions, which can play a tactical role, as riders often seek to escape this risk by moving to the front or breaking away entirely. Breakaway riders, however, encounter different risk factors, such as fatigue-induced impairments, and exposure to crosswinds that may affect stability.

In this paper, we aim to uncover the important role crashes play in deciding racing outcomes, and how competitive cyclists can make informed strategic decisions. We present a comprehensive mathematical approach that integrates energy consumption, drag reduction, and crash probability to determine the optimal position and timing for initiating a breakaway attempt. 

We begin in \S\ref{sec:governing_equations} by presenting the governing equations for our model. These comprise equations of motion for the riders and equations that capture the crashing probability. The equation of motion for each rider are based on Newton's second law, incorporating aerodynamic drag and power generation and losses due to rolling resistance and drivetrain friction. As we are interested in the role of attacking and escaping from the peloton, we incorporate a position-dependent drag coefficient, which depends on how far the rider is from the head of the peloton when riding in a group, and takes the same value as a rider at the front of the peloton when the rider is cycling alone. We introduce an objective function that balances the time advantage of winning a race with the probability of crashing.

To facilitate analytical progress, we choose to write the system in terms of the position of just two riders: the rider at the front of the peloton, and a special rider who will breakaway from the peloton at some point during the race. 

We first consider a simple case of a flat course with no crashes in \S\ref{sec:simple_attack}, to demonstrate the methodology. We extend this by incorporating the probability of a crash, in \S\ref{section:A simple breakaway with crash probability}. Here, we show how we may obtain an analytic expression for our objective function, which allows for easy optimization.

We develop the model further in \S\ref{sec:realistic_attack}, to account for rider fatigue. We also briefly discuss more realistic course profiles with variable elevation in \S\ref{sec:mountains}, and outline the numerical techniques required to analyse these more complex equations. In \S\ref{sec:discussion}, we draw conclusions from the analysis conducted. We use this to make relevant predictions to real-world cycle races, and we discuss extensions that would provide valuable further input for the optimization procedure.


\section{Governing equations}\label{sec:governing_equations}

\subsection{Dimensional equations}
We suppose that the race comprises a total of $N$ riders, and label these as $i=[1,N]$. We assume that all $N$ riders start off as a group, or \emph{peloton}, but that at some point during the race, a single rider can make a breakaway attempt. We assume that the riders are arranged in the peloton as a rectangular grid of cyclists, with $N_{\parallel}$ rows and $N_{\perp}$ columns, where $N=N_{\parallel}N_\perp$, and ordered such that rider $i=1$ is at the front of the peloton and rider $i=N$ is at the rear. 

We track the distance travelled along the course by each individual rider $i$ at time $\hat{t}$, and denote this by $\hat{x}_i(\hat{t})$. We begin by assuming that the course is flat (but generalize our approach for non-flat courses in Appendix~\ref{sec:mountains}). A force balance on an individual rider $i$ on flat terrain gives 
\begin{equation}\label{eq:individual_rider_dimensional}
\begin{split}
    \hat{m}_i \frac{\rd^2\hat{x}_i}{\rd \hat{t}^2}=\hat{P}_i(\hat{t})\left(\frac{\rd \hat{x}_i}{\rd \hat{t}}\right)^{-1}- \frac{1}{2}\hat{C}_{d,i}\rho A \left(\frac{\rd \hat{x}_i}{\rd \hat{t}}\right)^2.
    \end{split}
\end{equation}
In \eqref{eq:individual_rider_dimensional}, $\hat{m}_i$ is the mass of the rider and $\hat{P}(\hat{t})$ is the power the rider is exerting, which also accounts for rolling resistance and losses in the drivetrain. We choose to express the problem in terms of power exerted since most competitive cyclists carry a power-meter device, which they monitor closely during a race.  The constants $\rho$ and $A$ represent the air density and the rider's cross-sectional area, respectively. Finally, the drag coefficient $\hat{C}_{d,i}$ varies as a function of the distance between the individual rider and the first rider in the peloton. We model this functional dependence following Blocken \emph{et al.}~\cite{blockenAerodynamicDragCycling2018}, with drag decaying exponentially as we move further back into the peloton:
\begin{equation}\label{eq:dimensional_drag}
    \hat{C}_{d,i} = \begin{cases}
        \hat{C}_d^{\minm} + (\hat{C}_d^{\maxm}- \hat{C}_d^{\minm}) e^{-\lambda\frac{\hat{\zeta}_i}{d}} & \hat{\zeta}_i >0 \\
        \hat{C}_d^{\maxm} & \hat{\zeta}_i <0.
    \end{cases}
\end{equation}
Here, $d\approx 4 $ m is the axle-to-axle distance in a packed peloton, so that by definition $i=1+\hat{\zeta}_i/d$ is the drafting position of the individual rider within the peloton; 
$\lambda$ is a parameter that indicates how the drag force experienced by cyclists decays as they move into the back of the peloton. We estimate $\lambda\approx 0.25$ from the work of Blocken \emph{et al.}~\cite{blockenAerodynamicDragCycling2018}, and will assume this value throughout this paper. We also take a representative value of $\hat{C}_d^{\maxm}=0.9$ and $\hat{C}_d^{\minm}=0.05$ \cite{blockenAerodynamicDragCycling2018,MunsonBruceRoy1998Fofm}. 
We note that we have assumed that riders share identical aerodynamic properties, specifically the same cross-sectional area, $A$ and drag coefficient $C_{d,i}$, although the model readily generalizes to account for rider dependence of these parameters. 

Since the riders in the peloton all travel at the same speed, it is advantageous for us to work with the peloton behaviour on average. To study this, we introduce the definition of a peloton average of a quantity $f_i$,
\begin{align}
    \langle f \rangle=\frac{1}{N}\sum_{i=1}^N f_i,
\end{align}
Averaging \eqref{eq:individual_rider_dimensional} for the $N$ riders in the peloton then gives  
\begin{equation}\label{eq:peloton_dimensional}
\begin{split}
    \langle m\rangle \frac{\rd^2 \hat{x}_p}{\rd t^2}=\langle \hat{P}(\hat{t})\rangle  \left(\frac{\rd \hat{x}_p}{\rd \hat{t}}\right)^{-1} - \frac{1}{2}\langle \hat{C}_{d}\rangle &\rho A \left(\frac{\rd \hat{x}_p}{\rd \hat{t}}\right)^2 ,
    \end{split}
\end{equation}
where $\hat{x}_p=\langle \hat{x}\rangle+(N_\parallel-1)d/2$ defines the position of the front of the peloton. 
To close the model, we supply the initial condition $\hat{x}_p(0)=0$, which assumes that the race begins in peloton format, with the head of the peloton at the start line.

It is also useful to introduce the notion of energy depletion. If rider $i$ is exerting a power $\hat{P}_i(t)$, then by time $\hat{t}$ this rider will have consumed the following energy:
\begin{equation}
    \hat{E}_i(\hat{t})=\int_0^{\hat{t}} \hat{P}_i(\hat{t}')\,\rd\hat{t}'.
\end{equation}

From here on, we cease to consider $N$ equations for the position of all $N$ riders in the peloton, and instead focus on only two positions: the head of the peloton, $\hat{x}_p(\hat{t})$, and one `special' rider, whose race behaviour we wish to optimize. Without loss of generality, we label the special rider by index $i$. Their position is then given by $\hat{x}_i(\hat{t})$, and $\hat{x}_i(0)=-\hat{\zeta}_i$, so that this rider starts behind the head of the peloton.

\subsection{Non-dimensionalization}\label{sec:scalings}
To simplify the analysis of the equations of motion we introduce the following dimensionless variables: 
\begin{align}
\hat{x}=Lx, \qquad \hat{\zeta}_i=d\zeta_i, \qquad   \hat{t}=\left(\frac{\langle \hat{C}_d\rangle \rho AL^3}{2\langle\hat{P}(0)\rangle}\right)^{1/3}t, \nonumber \\
    \hat{P}=\langle\hat{P}(0)\rangle P, \hspace{7mm}
 \hat{E}=\left(\frac{\langle \hat{C}_d\rangle \rho A \langle\hat{P}(0)\rangle L^3}{2}\right)^{1/3}E.
\end{align}
With this choice of scalings, the course length is one unit and the course will be completed in one time unit by the rider at the front of the peloton. \jcvc{Moreover, the average energy expenditure of the peloton will also be equal to one.} 

The dimensionless counterpart of the equation of motion for the peloton, \eqref{eq:peloton_dimensional}, is then given by 
\begin{equation}\label{eq:dimensionless_peloton}
    \epsilon \frac{\rd^2 x_p }{\rd t^2}=P(t)\left( \frac{\mathrm{d}x_p}{\mathrm{d}t}\right)^{-1}-\left(\frac{\mathrm{d}x_p}{\mathrm{d}t}\right)^2,
\end{equation}
where 
\begin{align}
\epsilon & = \frac{2\langle m\rangle}{L \langle C_d\rangle \rho A}\approx 0.005.
\end{align}
Here, we have used $\rho=1.225$ kg$\cdot$m$^{-3}$, $\langle m\rangle=70$ kg, $L\approx 100$ km and $\langle C_d\rangle A \approx 0.4$ m$^2$ \cite{vandruenenAerodynamicAnalysisUphill2021}.
\jcvc{We may interpret the smallness in $\epsilon$ by noting that changes in velocity happen very quickly when compared to the duration of the race.}
The dimensionless counterpart to the equation of motion for the individual rider, \eqref{eq:individual_rider_dimensional}, is
\begin{equation}\label{eq:dimensionless_individual_rider}
    \epsilon m_i \frac{\rd^2 x_i}{\rd t^2}={P_i(t)}\left(\frac{\rd{x}_i}{\rd t}\right)^{-1}-C_{d,i}\left(\frac{\rd x_i}{\rd t}\right)^2,
\end{equation}
where $m_i={\hat{m}_i}/{\langle m\rangle}$ is the dimensionless mass of the rider, and $C_{d,i}=\hat{C}_{d,i}/\langle \hat{C}_d\rangle$ is the dimensionless drag coefficient. \jcvc{For completeness, we list the initial conditions for the peloton and the individual rider:}
\begin{equation}
    x_p(0)=0,\quad x_i(0)=-\delta\zeta_i,\quad \frac{\rd{x}_p(0)}{\rd t}=\frac{\rd{x}_i(0)}{\rd t}=1.
\end{equation}
where $\delta=d/L \ll 1$. For a typical course length of 150km, $\delta \approx 2\times 10^{-5}$. In all subsequent analysis, we will neglect terms of order $\delta$, except in Appendix \S\ref{sec:transition_layers}, where we analyse the behaviour during an escape from the peloton, and must consider the system on the peloton lengthscale.    

The dimensionless energy is given by 
\begin{align}
\label{eq:dimensionless energy}
E(t)=\int_0^tP_i(t')\,\rd t',
\end{align}
where ${\tf}$ is the time it takes for the individual rider to complete the course. We will pay particular attention to the dimensionless quantity $E^\ast=E(t_f)$, which measures the total energy used during the entire race by the special rider compared with the peloton average. 

Motivated by the small value of the inertia, $\epsilon\ll1$, we will focus on the quasi-steady regime where $\epsilon=0$. In physical terms, this means that changes in power $P_i(t)$ drive instantaneous jumps in velocity. Whilst this is not true when looking at the microscale, it is an excellent model when comparing the time it takes to accelerate/decelerate to a new velocity (seconds) with the time it takes to complete the entire course (hours). The microstructure associated with inertial effects and acceleration is explored further using boundary layer methods in Appendix \S\ref{sec:transition_layers}.
\subsection{Crashing}

Cycling near the rear of the peloton has clear advantages, chiefly a drag reduction close to 90\% \cite{blockenAerodynamicDragCycling2018}. However, it is also associated with risks. In particular, if there is is a crash at a given location in the peloton, since a crash propagates backwards in the peloton, a rider who is further back is more likely to crash. Therefore, we need to model the probability of crashing as part of our race-winning strategy. To find the risk penalty associated with crashing, we use a conditional-probability approach. 

\begin{figure}[t]
\centering
\includegraphics[width=0.99 \linewidth]{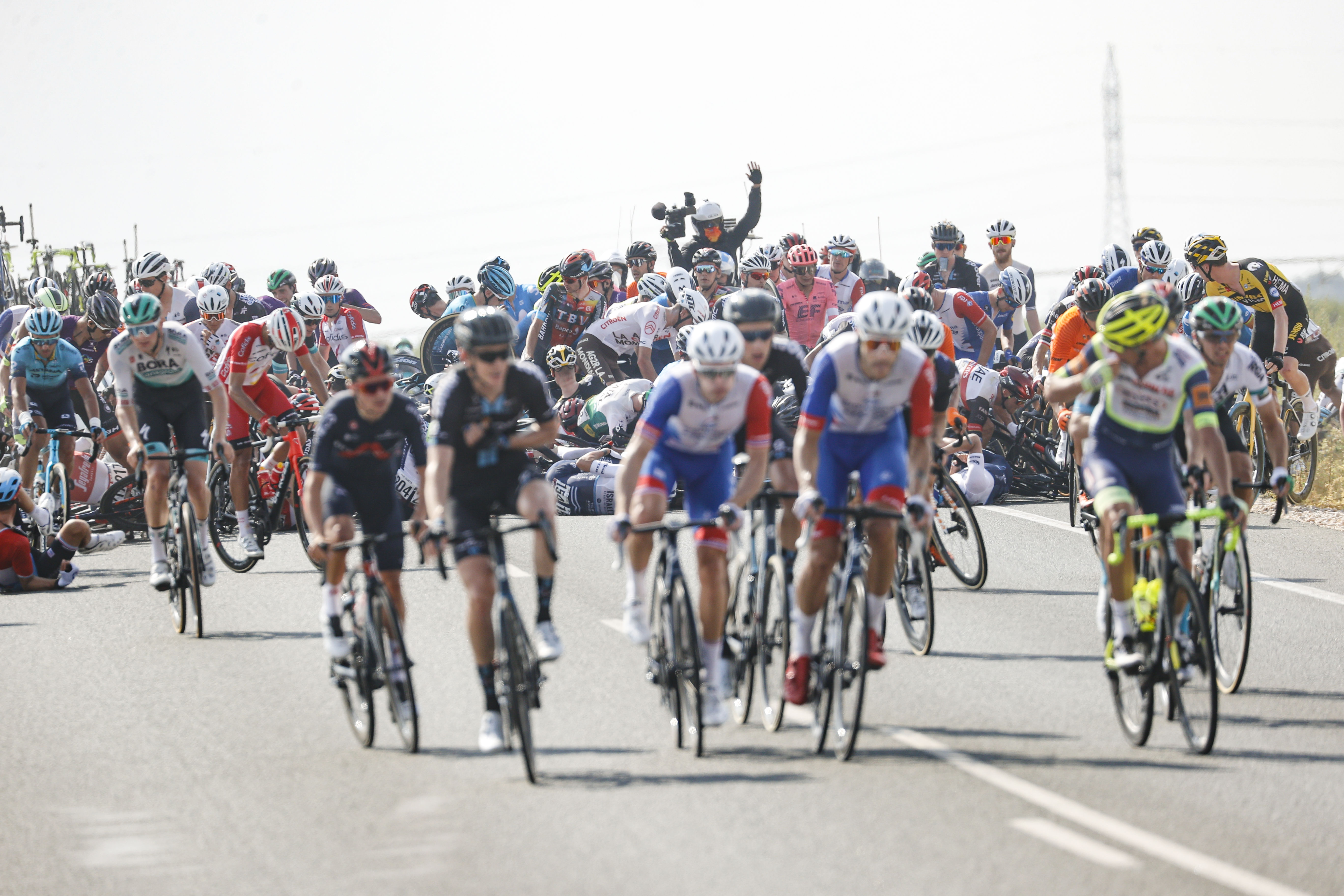}
\caption{The immediate aftermath of a crash (\textit{pile-up}). Riders positioned ahead of the start of the crash are spared. The crash propagates backwards through the peloton. Credit: Luis \'{A}ngel G\'{o}mez.}
\label{fig:crash}
\end{figure}

\jcvc{Let $C(x)$ denote the event of a crash taking place between $x$ and $x+\rd x$. 
We define $C_i\vert C(x)$ as the event that the rider in position $i$ crashes given that a crash occurred at location $x$. In other words, given that a crash has taken place at $x$, $C_i\vert C(x)$ is the event that rider $i$ is involved. Moreover, we introduce $S_k\vert C(x)$ as the event that the first rider to crash is at position $k$, given that a crash occurred at location $x$. In order to compute the probability $\mathbb{P}(C_i(x))$, we consider $\mathbb{P}(C_i\vert S_k(x))$, the probability that rider $i$ is involved in the crash given that rider $j$ is involved in it.} We assume that any riders ahead of the crash are unaffected, but that riders behind may crash, with a probability that decays exponentially as we move away from the crash: 
\begin{equation}\label{eq:CrashPropagation}
    \mathbb{P}((C_i\vert C)\vert (S_k\vert C))=\mathbb{P}(C_i\vert S_k)=\begin{cases}0 &\text{if }i<
    k,\\
    e^{-\omega(i-k)}& \text{if }i\geq k.
        
    \end{cases}
\end{equation}
Here, $\omega(x)$ is a parameter that measures how far back a crash propagates on average, and in general will depend on the peloton speed, road conditions (for example, rain and terrain type), along with other factors. Generally, we will set $\omega<{1}$, so that an average crash will propagate over $\omega^{-1}>1$, so more than one rider. We also note that $\mathbb{P}(C_i\vert S_i)=1$, as required. 
We assume that the crash is equally likely to start at any position, so that $\mathbb{P}(S_k\vert C)=1/N$. Non-uniform distributions for $\mathbb{P}(S_k\vert C)$, as well as a general model for the crash propagation \eqref{eq:CrashPropagation} are explored in \S\ref{sec:complex_crash}. Invoking the law of total probability, we obtain the probability that the rider at position $i$ crashes given a crash has occurred:
\begin{align}
   \mathbb{P}(C_i|C(x))&=\sum_{k}\mathbb{P}((C_i\vert C)\vert (S_k\vert C))\mathbb{P}(S_k\vert C) \nonumber  \\ &=\frac{1}{N}\sum_{k=1}^je^{-\omega (i-k)}=\frac{1}{N}\frac{1-e^{-\omega i}}{1-e^{-\omega}}\equiv H(i;\omega).
   \label{eq:crash probability given someone crashes}
\end{align}
The probability of rider $i$ being involved in a crash at some point in the race (defined as $\mathcal{P}$) follows from a second application of the law of total probability,
\begin{equation}\label{eq:ProbDefn}
    \mathcal{P}\equiv\mathbb{P}(C_i)=\int_0^1  \mathbb{P}(C_i|C(x))\mathbb{P}(C(x))\,\rd x,
\end{equation}
where $\mathbb{P}(C(x))$ is the probability of a crash occurring between $x$ and $x+\rd x$. For simplicity, here we will assume that the probability of crashing is independent of $x$ (no cobblestones, dangerous corners or other treacherous features), so that 
\begin{equation}\label{eq:ProbArbitrary}
    \mathcal{P} = \mathbb{P}(C)\int_0^1H(i(x);\omega)\,\rd x,
\end{equation}
where we write $i(x)$ to explicitly denote that the position of the rider within the peloton changes along the course of the race. 
To estimate the crash density $\mathbb{P}(C(\hat{x}))$, we consider a 21-stage, 3500 km-long Grand Tour. There are typically 1--3 crashes per stage, depending on environmental conditions and course difficulty~\cite{decock_incidence_2016}. This corresponds to a crash every 50--100 km of racing, and so $\mathbb{P}(C(\hat{x}))\approx 0.01-0.02$ crashes/km. In dimensionless variables, $\mathbb{P}(C(x))=L\mathbb{P}(C(\hat{x})=1-2$ crashes/unit dimensionless distance, which is the same as $\mathbb{P}(C(x))=1-2$ crashes/stage. Thus, we use $\mathbb{P}(C)=2$ unless otherwise stated. 

\subsubsection{A simple example}\label{sec:CrashingSimpleExample}
Here, we consider an example computation of $\mathcal{P}$ that will be useful for later computations. In particular, in \S\ref{sec:simple_attack} we will consider the simple scenario with the rider initially in the peloton, until they break away at attack position $x=x_a$  (so that $i(x)=i$ for $x<x_a$, and $i(x)=1$ for $x>x_a$). The integral in \eqref{eq:ProbArbitrary} thus evaluates to:
\begin{equation}\label{eq:CrashingProbabilitySimpleAttack}
\begin{split}
    \mathcal{P}&=\mathbb{P}(C)\left[x_a H(i;\omega)+(1-x_a)H(1;\omega)\right]\\
    &=\frac{\mathbb{P}(C)}{N}\left(x_a\frac{1-e^{-\omega i}}{1-e^{-\omega}}+1-x_a\right).
    \end{split}
\end{equation}

\subsection{Objective function  and optimization}
The rider wishes both to finish ahead of the peloton and to avoid crashing. Different riders will be willing to accept more or less risk. To this end, we introduce an objective function to evaluate the suitability of a racing strategy. The objective function should be low when both the rider finishes far ahead of the peloton and their probability of crashing is small. These considerations motivate the definition of the objective function $\mathcal{M}\in\mathbb{R}$ as:
\begin{align}\label{eq:metric_minima_definition}
    \mathcal{M} = -\beta \Deltat +(1-\beta)\mathcal{P},
\end{align}
where $\Deltat={\tf}-{\tp}$, ${\tf}$ is the time that the special rider $i$ takes to finish, ${\tp}$ is the time it takes for the peloton to complete the course, and $\mathcal{P}$ is the crashing probability defined in \eqref{eq:ProbDefn}. Under our choice of scalings, ${\tp}=1$. The parameter $\beta$ is a measure of the risk the special rider is willing to assume: $\beta=0$ corresponds to risk-averse riders, who prioritize not crashing over a larger winning margin, while $\beta=1$ applies to riders who are willing to take risks at any cost to maximize the win margin. Inclusion of a risk parameter is crucial to the applicability of the model. For example, the \textit{Maillot Jaune} in the later stages of a Grand Tour should choose a low value of $\beta$, as they are already winning the race and have much to lose if they crash and are forced to retire.

The goal of this paper is to find the strategy for rider $i$ that minimizes $\mathcal{M}$ for a given risk tolerance $\beta$. By strategy, we mean the allocation of resources from the energy budget for cyclist $i$, $E^*$, by changing the power $P_i(t)$ during the race \cite{foster_regulation_2005}. 



\section{A simple breakaway with no crashes}\label{sec:simple_attack}
\begin{figure*}
    \centering
    \includegraphics[width=0.99\linewidth]{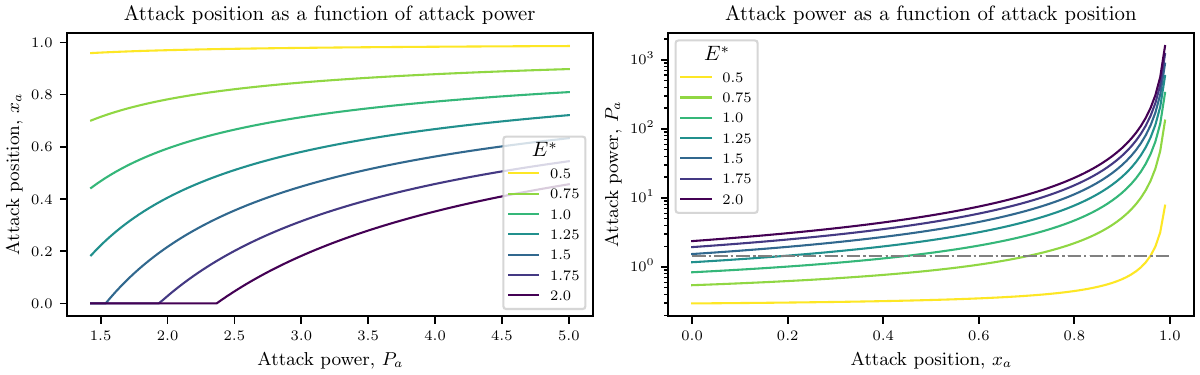}
    \caption{Left: plots of the attack position as a function of attack power for different energy budgets $E^*$, given by \eqref{eq:xa in terms of Pa and E*}. Right: plots of the attack power as a function of the attack position for different energy budgets $E^\ast$, given by \eqref{eq:Pa in terms of xa and E}. The dashed grey line represents the minimum power required for the attack to be successful, $P_a^\minm=C_{d,1}$. For both plots the initial position of the rider inside the peloton is $i=5$, which translates to $P_{l,i}=C_{d,i}=0.46$. Other parameters are $\omega=0.5$ and $C_{d,1}=1.43$.  }    
    \label{fig:x_a_P_a_plots}
\end{figure*}

To illustrate the mathematical model, we begin in this section by considering a race in which no crashes take place. We will suppose that a single rider sits at position $i$ in the peloton until making a breakaway attempt at some attack point along the course, which we denote as $x_a$. By our choice of non-dimensionalization, the peloton moves at unit dimensionless speed, and hence $x_p(t)=t$. Second, we assume the power exerted has the following form. Before the attack, the rider is in the peloton at position $i$ and so the power exerted is $P_i=C_{d,i}$, using \eqref{eq:dimensionless_individual_rider}. We denote $P_{l,i}=C_{d,i}$ as the lurking power of the rider at position $i$ in the peloton. We note that, due to the choice in non-dimensionalization, $P_{i,l}<1$ for $i>1$. After the rider makes their breakaway move at, say, $t={\ta}$, we suppose that they provide a constant power $P_a>P_{l,i}$ until the end of the race. The total energy expenditure is then simply 
\begin{equation}\label{eq:energy_computation_flat_simple}
E^\ast=P_{l,i}{\ta}+P_a({\tf}-{\ta}).
\end{equation}

As the peloton moves with unit speed, the attack position is related to the attack time by 
$x_a = {\ta}$ (neglecting terms in $\zeta_i$ as noted earlier). After the attack, the time left to complete the course is 
${\tf}-{\ta}=({1-x_a})/{v_a}$ where $v_a \equiv \rd x_i/\rd t = \left(P_a/C_{d,1}\right)^{1/3}$ is the attack velocity, using \eqref{eq:dimensionless_individual_rider} and the fact that the drag on a solo rider is $C_{d,1}$. We can thus find the earliest attack location for a given attack power, $P_a$, that a rider can provide and their total energy reserve, $E^\ast$, 
\begin{align}
\label{eq:xa in terms of Pa and E*}
    x_{a}=\max\left\{\frac{C_{d,1}^{1/3}P_a^{2/3}-E^\ast}{C_{d,1}^{1/3}P_a^{2/3}-C_{d,i}},0\right\}. 
\end{align}
Here, we have introduced the $\max$ operator to ensure the attack position is in the domain, i.e. $x_a>0$. Since we require $v_a>1$ for an attack to be successful, this gives a minimum attack power, $P_a^\minm=C_{d,1}$ and thus an associated earliest attack time, 
\begin{align}\label{eq:minimum_attack_position}
    x^{\mathrm{min}}_{a}=\max\left\{\frac{C_{d,1}-E^\ast}{C_{d,1}-C_{d,i}},0\right\},
\end{align}
using \eqref{eq:xa in terms of Pa and E*}. For the rider to have a chance of winning, we must have $x^{\mathrm{min}}_{a}<1$, and so this means that the minimum energy of the rider for an attack to be successful is $E^\ast_{\mathrm{min}}=C_{d,i}$.

\jcvc{In the left panel of fig.~\ref{fig:x_a_P_a_plots} we plot the dependence of $x_a$ on the attack power for different values of $E^*$, given by \eqref{eq:xa in terms of Pa and E*}. For small energy budgets, the earliest attack position is close to the finish line, but for higher energy budgets the attack position is closer to the start line. In particular, if the energy budget is high enough (namely if $E^*>C_{d,1}$), then the earliest attack position $x_a\leq0$. In this case, this means that a rider has sufficient energy to make a break at the beginning of the race (and still have surplus energy at the end of the race)}.

Alternatively, a racer may have identified a particular place where they wish to breakaway. In this case, we can rearrange \eqref{eq:xa in terms of Pa and E*} for $P_a$ to determine the required power for the optimal breakaway that uses all their energy reserves by the end of the race:
\begin{align}\label{eq:Pa in terms of xa and E}
    P_a=\lr{\frac{E^\ast-C_{d,i}x_a}{C_{d,1}^{1/3}(1-x_a)}}^{3/2}.
\end{align}
\jcvc{The dependence of $P_a$ on the attack position $x_a$ is plotted in the right panel of fig.~\ref{fig:x_a_P_a_plots}. We show the minimum attack power required to escape from the peloton ($P_a>C_{d,1}$ from $v_a>1$), as a dashed grey line.}

\begin{figure*}
    \centering
    \includegraphics[width=0.99\linewidth]{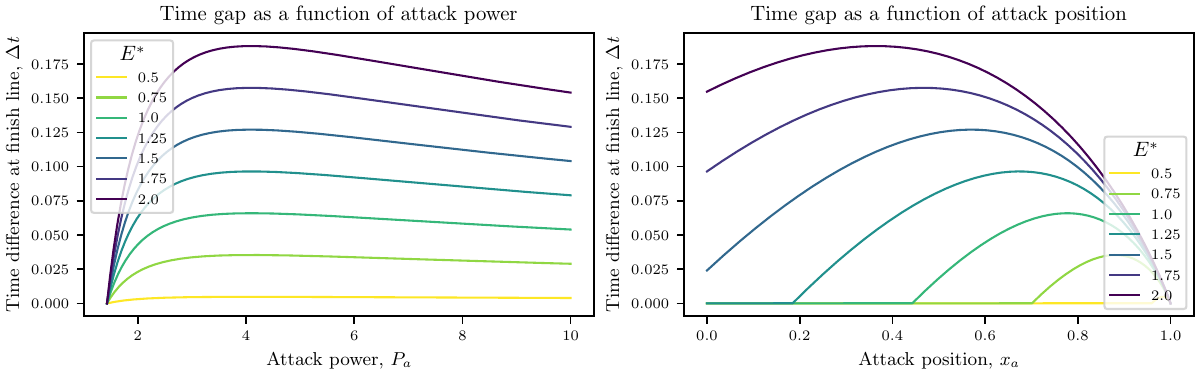}
    \caption{Dependence of finish time of special rider ahead of the peloton for a fixed energy budget $E^\ast$, as a function of either attack power, $P_a$, given by \eqref{eq:Delta T for flat course in terms of Pa} (left) or attack position, $x_a$, given by \eqref{eq:Delta T for flat course in terms of xa} (right). For both plots the initial position of the rider inside of the peloton is $i=5$, which translates to $P_{l,i}=C_{d,i}=0.46$. Other parameters are $\omega=0.5$ and $C_{d,1}=1.43$.   }
    \label{fig:Delta_T_plots}
\end{figure*}

We may also solve for the finish time ahead of the peloton for the breakaway rider in terms of $P_a$ and $E^\ast$, 
\begin{align}
\label{eq:Delta T for flat course in terms of Pa}
\Deltat=1-{\tf}=\frac{C_{d,i}-E^\ast}{C_{d,1}^{1/3}P_a^{2/3}-C_{d,i}}\lr{\lr{\frac{C_{d,1}}{P_a}}^{1/3}-1},
\end{align}
or in terms of $x_a$ and $E^\ast$,
\begin{align}\label{eq:Delta T for flat course in terms of xa}
    \Deltat=1-x_a-\frac{(1-x_a)^{3/2}C_{d,1}^{1/2}}{(E^\ast-C_{d,i}x_a)^{1/2}}.
\end{align}
We may see how the finish time ahead of the peloton, $\Deltat$ varies with either the attack power or the breakaway position in fig.~\ref{fig:Delta_T_plots}. In the left panel, we observe the that the time gap initially increases as the breakaway power rises, but eventually it reaches a maximum before falling again.  Similarly, the time gap initially grows with attack position before reaching a maximum before falling as the attack position gets closer to the end of the course with $\Deltat=0$ at $x_a=1$. Another feature we highlight from the plots for $\Deltat(x_a)$ is the discontinuous derivative of the function at the minimum attack position. Breakaways attempted before this position lack the power necessary to overcome the increased drag when riding in isolation, and the rider simply returns to the peloton and finishes the race with the main group, with $\Deltat =0$. 

A key observation from fig.~\ref{fig:Delta_T_plots} is the existence of an optimum that maximizes the time difference between the breakaway rider and the peloton for a given attack power $P_a$ or attack position $x_a$, when taking account of the available energy reserves. In the following section, we will study the dependence of this optimum on the parameters, and study how this is affected by the introduction of a non-zero crashing probability. 

\section{A simple breakaway with crash probability}
\label{section:A simple breakaway with crash probability}

\subsection{Assembling the objective function}
We now extend the model to study the effect that a non-zero probability of crashing has on the existence of a strategy that maximizes the finish time ahead of the peloton. %

Substituting for $\Deltat$ and $\mathcal{P}$ in \eqref{eq:metric_minima_definition} using \eqref{eq:Delta T for flat course in terms of xa} and \eqref{eq:CrashingProbabilitySimpleAttack}, respectively, gives 
\begin{equation}\label{eq:MetricSimpleAttack}
\begin{split}
    \mathcal{M}&=-\beta\left(1-x_a-\frac{(1-x_a)^{3/2}C_{d,1}^{1/2}}{(E^\ast-C_{d,i}x_a)^{1/2}}\right)\\ &+(1-\beta)\frac{\mathbb{P}(C)}{N}\left(x_a\frac{1-e^{-\omega i}}{1-e^{-\omega}}+1-x_a\right).
    \end{split}
\end{equation}

\begin{figure*}
    \includegraphics[width=0.99\textwidth]{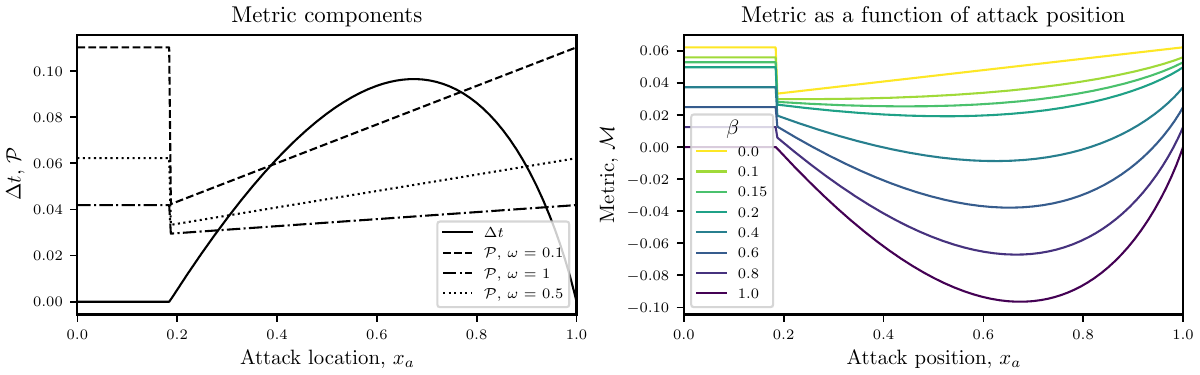}
    \caption{Left: the finish time ahead of the peloton, $\Deltat$, and the crashing probability $\mathcal{P}$ plotted as a function of the attack position $x_a$, given, respectively, by  \eqref{eq:Delta T for flat course in terms of xa} and \eqref{eq:CrashingProbabilitySimpleAttack}. We plot $\mathcal{P}$ for $\omega=0.1, 0.5$ and $1$. Right: objective function $\mathcal{M}(x_a)$, given by \eqref{eq:metric_minima_definition}, plotted as a function of attack position, $x_a$. Parameter values are $E^\ast=1.2$, $i=5$, $C_{d,i}=0.46$, $C_{d,1}=1.43$, $\omega=0.5$,  $\mathbb{P}(C)=2$ and $N=75$. 
    }
    \label{fig:metric_example}
\end{figure*}

In the left panel of fig.~\ref{fig:metric_example}, we plot the components that make up the objective function, $\Deltat$ and $\mathcal{P}$, given, respectively, by \eqref{eq:Delta T for flat course in terms of xa} and \eqref{eq:CrashingProbabilitySimpleAttack}, as a function of the breakaway point, $x_a$. In the right panel of fig.~\ref{fig:metric_example}, we plot the objective function $\mathcal{M}$ as a function of the attack position for several values of the risk index level $\beta$. 

For suitably large values of $\beta$, the value of $x_a$ that minimizes $\mathcal{M}$ is a local minimum. However, when $\beta$ is smaller, the minimum value of $\mathcal{M}$ occurs at the discontinuity in the graph, which occurs at $x_a=x_a^\minm$, defined by \eqref{eq:xa in terms of Pa and E*}. In the next section, we will analyse this further, to find a method to determine the value of $x_a$ where $\mathcal{M}$ is minimized, and hence the optimal attack position.

\subsection{Results}\label{sec:SimpleAttackResults}
We now proceed to systemically analyse the dependence of the critical attack position, and the finish time ahead of the peloton associated with such a location, on the system parameters. As seen earlier, the minimum value of $\mathcal{M}$ can either occur at a local minimum in $\mathcal{M}(x_a)$ or at a discontinuity in $\mathcal{M}(x_a)$. 

If $\mathcal{M}$ attains its minimum value at a local minimum, say $x_a^\dagger\in (0,1)$, we may find the optimal attack position that minimizes the objective function by solving $\partial\mathcal{M}/\partial x_a=0$,
which may be written as 
the depressed cubic 
\begin{equation}\label{eq:DepressedCubic}
\begin{split}
    &(1-\beta)\left[N H(i;\omega)-1\right]\frac{\mathbb{P}(C)}{N}+\beta \\&\hspace{30mm}+\beta C_{d,1}^{1/2}\left[\frac{C_{d,i}}{2}\eta^3-\frac{3}{2}\eta\right]=0,
    \end{split}
\end{equation}
where 
\begin{equation}\label{eq:EtaDefn}
    \eta=\left(\frac{1-x_a^\dagger}{E^\ast-C_{d,i}x_a^\dagger}\right)^{1/2}.
\end{equation}
While \eqref{eq:DepressedCubic} may be solved analytically for $\eta$ and hence $x^\dagger_a$, this is not required. Instead, we invert \eqref{eq:EtaDefn} for $x_a^\dagger$, to give
\begin{equation}
\label{eq:xa in terms of eta}
    x_a^\dagger=\frac{1-E^\ast \eta^2}{1-C_{d,i}\eta^2},
\end{equation}
and since the solution for $\eta$ from \eqref{eq:DepressedCubic} is independent of $E^\ast$, we immediately see the optimal attack position changes linearly with the budget energy $E^\ast$. We may also find the corresponding optimal attack power, say $P_a^\dagger$, in terms of $\eta$ by making use of \eqref{eq:Pa in terms of xa and E},
\begin{equation}\label{eq:OptimalAttackPower in terms of eta}
    P_a^\dagger = \frac{1}{C_{d,1}^{1/2}\eta^3}.
\end{equation}
This optimal attack power is thus independent of the energy budget. 

The minimum of $\mathcal{M}$ may also be located at the discontinuity at $x_a=x_a^\minm$, defined by \eqref{eq:xa in terms of Pa and E*}. To find the global minimum of  $\mathcal{M}$, we thus compare the interior minimum $x_a^\dagger$ found from \eqref{eq:xa in terms of eta} with the value of $\mathcal{M}$ at $x_a=x_a^\minm$ given by \eqref{eq:xa in terms of Pa and E*}. The optimal location, $x_a^\ast$, is then given by
\begin{equation}
    x_a^\ast=\mathrm{argmin}\left\{\mathcal{M}(x_a^\dagger),\mathcal{M}(x_a^\minm)\right\}.
\end{equation}
Similarly, once the optimal attack location is known, we may use \eqref{eq:Pa in terms of xa and E}, to find the optimal attack power, $P_a^\ast$. We now study how the optimal attack location changes as the risk index is varied.

\subsubsection{Role of risk}
\begin{figure*}
    \includegraphics[width=0.99\textwidth]{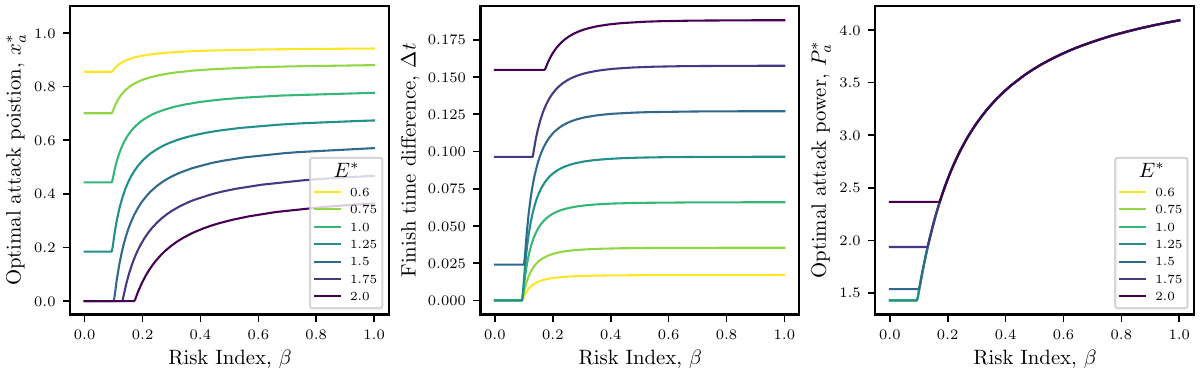}
    \caption{Dependence on assumed risk. Left: attack position $x_a^\ast$ that minimizes the objective function as a function of the risk index $\beta$. Middle: Finish time ahead of the peloton when the break is attempted at the optimal position as a function of the risk index. Right: Optimal attack power computed from $x_a^\ast$ using \eqref{eq:Pa in terms of xa and E}. Optimal solutions are shown for several values of the energy budget $E^\ast$. Parameter values are $i=5$, $C_{d,i}=0.46$, $C_{d,1}=1.43$, $\lambda=0.25$, $\omega=0.5$, $\mathbb{P}(C)=2$ and $N=75$.
    }
    \label{fig:x_optimal_vs_risk}
\end{figure*}

In the left panel of fig.~\ref{fig:x_optimal_vs_risk} we plot how the optimal attack position, $x_a^\ast$, varies as a function of the assumed risk $\beta$, for several values of the energy budget. We observe the aforementioned jump from $x_a^\minm$ to $x_a^\dagger>0$ as $\beta$ is increased. We see that $x_a^\ast=0$ when $\beta=0$ if $E^\ast >E^\ast_{\mathrm{crit}}=C_{d,1}$, using \eqref{eq:minimum_attack_position}, and $x_a^\ast>0$ when $\beta=0$ if $E^\ast <E^\ast_{\mathrm{crit}}=C_{d,1}$. 
For the parameters used in fig.~\ref{fig:x_optimal_vs_risk}, we find $E^\ast_{\mathrm{crit}}=1.43$. we also see that the optimal attack position asymptotes to an $E^\ast$-dependent limit value as $\beta\rightarrow1$. 

In the middle panel of fig.~\ref{fig:x_optimal_vs_risk} we plot how the finish time ahead of the peloton varies as a function of risk, when the rider makes an attack at the optimal location $x_a^\ast$ (as given in the left panel). For sufficiently high energy budgets, if the rider breaks at the start of the race they will win, but for lower energy budgets, $E^\ast\lesssim E^\ast_{\mathrm{crit}}$, the rider does not possess sufficient energy to break from the start, and ends up finishing the race with the peloton. This highlights how, when the risk index $\beta$ is low, the best strategy is to break away from the peloton not to win but to avoid crashing. As the risk increases and the optimal attack position moves away from $x_a^\ast=0$, the finish time ahead of the peloton subsequently increases as well, converging to an $E^\ast$-dependent limit as $\beta\rightarrow1$. 

It is especially noteworthy that it is possible to win with an energy budget that is appreciably lower than that used on average in the peloton, namely $E^\ast\approx 0.6$, but to do so requires the rider to be willing to assume a more elevated level of risk (right panel of fig.~\ref{fig:x_optimal_vs_risk}). More generally, riders with an energy budget $E^\ast<E_{\mathrm{crit}}^\ast$ must assume a minimum level of risk, about $\beta\approx0.1$ to win the race.  This is consistent with the intuition that a strong competitor can assume a much lower risk to win, while a weaker opponent must be willing to race with a riskier strategy.

In the right panel of fig.~\ref{fig:x_optimal_vs_risk} we plot the optimal attack power $P_a^\ast$, computed from $x_a^\ast$ using \eqref{eq:Pa in terms of xa and E}. We observe that when $x^\ast_a>0$, the optimal power for different energy budgets collapses into a single curve, as predicted in \eqref{eq:OptimalAttackPower in terms of eta}, where we theorized that $P_a^\ast$ was independent of $E^\ast$ provided the optimal attack position was non-zero.  

\subsubsection{Dependence on energy budget}

Intuition tells us that the energy budget will play an important role on racing strategy. In particular, we found in \eqref{eq:xa in terms of eta} that the optimal attack position varies linearly with the energy budget, provided we are in a region of parameter space where the objective function is smooth. 

\begin{figure*}
    \centering
    \includegraphics[width=0.99\linewidth]{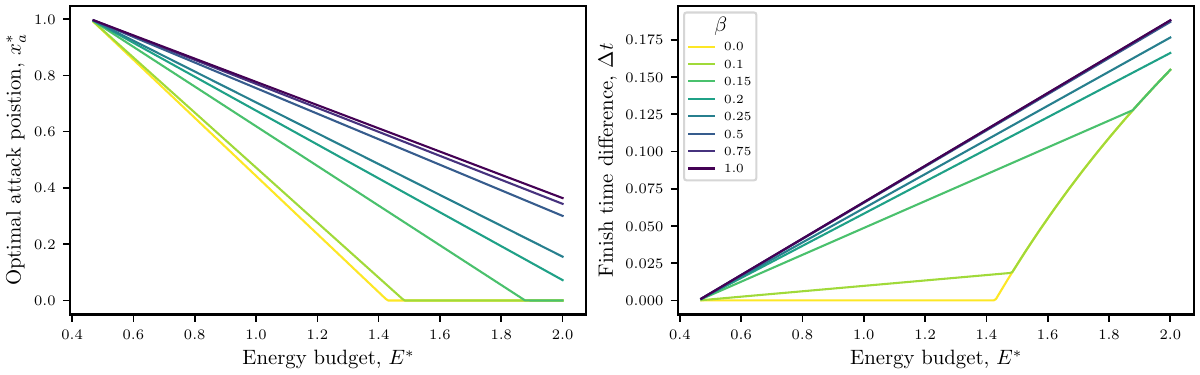}
    \caption{Optimal attack position (left) and associated finish time ahead of the peloton (right) plotted as a function of energy budget $E^\ast$ for several risk indices. Parameter values are $i=5$, $C_{d,i}=0.46$, $C_{d,1}=1.43$, $\lambda=0.25$, $\omega=0.5$, $\mathbb{P}(C)=2$ and $N=75$. }
    \label{fig:x_optimal_vs_energy}
\end{figure*}

In the left panel of fig.~\ref{fig:x_optimal_vs_energy}, we plot how the optimal attack position changes with $E^\ast$ for various risk levels. For low energy budgets the optimal location to attack is as close as possible to the end of the race. As the energy budget increases, the optimal location decreases linearly (as predicted by \eqref{eq:xa in terms of eta}) until it reaches $x_a=0$ again for conservative riders with high energy budgets, where the rider is strong enough to break from the start and win as well. Cyclists with high energy budgets who prioritize winning are recommended to attack close to the midpoint of the race. 

In the right panel of fig.~\ref{fig:x_optimal_vs_energy}, we plot how the finish time ahead of the peloton (when breaking at the optimal location) changes with the energy budget. As the energy budget increases, the time difference will also increase (as intuition suggests). For more conservative riders, there is a sharp transition at a critical value of $E^\ast$, where the rider can win by attacking at the start of the race.

\subsubsection{Minimum risk to win}
\begin{figure*}
    \centering
    \includegraphics[width=0.99\linewidth]{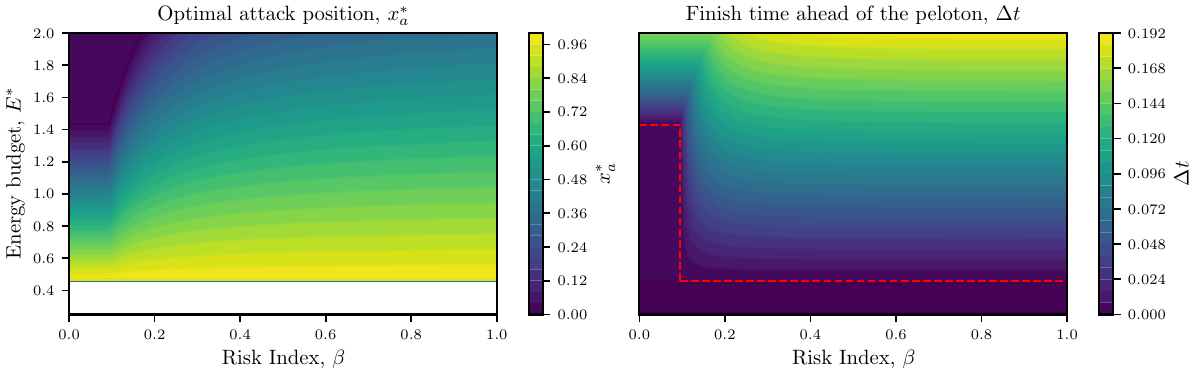}
    \caption{Left: contours of the optimal attack position as a function of the risk index $\beta$ and the energy budget $E^\ast$. For $E^\ast<E^\ast_\minm=C_{d,i}$ (in white) there is no attack position $0 \le x_a \le 1$ that leads to a race win. 
    Right: contours of the finish time ahead of the peloton $\Deltat(\beta, E^\ast)$ when attacking at the location given by right panel. The dashed red line is the function $\mathcal{E}_\minm(\beta)$ from \eqref{eq:CriticalRisk}. Parameter values are $i=5$, $C_{d,i}=0.46$, $C_{d,1}=1.43$, $\lambda=0.25$, $\omega=0.5$, $\mathbb{P}(C)=2$ and $N=75$. 
    }
    \label{fig:x_a_optimal_vs_energy_and_risk}
\end{figure*}

In fig.~\ref{fig:x_a_optimal_vs_energy_and_risk} we show contours of the optimal attack position (left) and the finish time ahead of the peloton (right) as a function of the risk index and the energy budget. 
From the right panel we can see how each risk index $\beta$ has an associated minimum energy budget that allows the breakaway rider to win the race, which we define as $\mathcal{E}_\minm(\beta)$. We see that $\mathcal{E}_\minm(\beta)$ is a piecewise constant function, 
\begin{align}
    \mathcal{E}_\minm(\beta)=\begin{cases}
        E^\ast_{\mathrm{crit}}=C_{d,1} & \beta <\beta^\ast,\\
    E^\ast_\minm = C_{d,i} & \beta > \beta^\ast.
    \end{cases}
\end{align}
This corresponds to the optimal strategy of either attacking at the beginning of the race, $x_a=0$ (when $\beta<\beta^\ast$) or attacking at some intermediate point in the race $x_a>0$ (when $\beta>\beta^\ast$). We overlay a plot $\mathcal{E}_\minm(\beta)$ in the right panel of fig.~\ref{fig:x_a_optimal_vs_energy_and_risk} (right).

We compute the critical risk index $\beta^\ast$ in terms of the system parameters in Appendix~\ref{sec:minimum_risk_appendix}, which gives
\begin{equation}\label{eq:CriticalRisk}
    \beta^\ast=\frac{\frac{\mathbb{P}(C)}{N}\left(\frac{1-e^{-\omega i}}{1-e^{-\omega}}-1\right)}{\frac{\mathbb{P}(C)}{N}\left(\frac{1-e^{-\omega i}}{1-e^{-\omega}}-1\right)+\frac{1}{2}\left(1-\frac{C_{d,i}}{C_{d,1}}\right)}.
\end{equation}
For the parameter values used to produce fig.~\ref{fig:x_a_optimal_vs_energy_and_risk}, the critical risk is $\beta^\ast=0.0949$. In practice, this means that riders who take greater risks require less energy to win, with a sharp cut-off at $\beta=\beta^\ast$. 

For an elite rider, it is perhaps more interesting to consider the converse question, that is: for a given energy budget, what is the smallest risk one must take to win the race? This is provided by the inverse function, $\beta_{\minm}(E^\ast)$, given by 
\begin{align}
    \beta_{\minm}(E^\ast)=\begin{cases}
        \beta^\ast & E^\ast_\minm < E^\ast<E^\ast_{\mathrm{crit}},\\
        0 & E>E^\ast_{\mathrm{crit}},
    \end{cases}
\end{align}
where $E^\ast_{\mathrm{crit}}=C_{d,1}$, $E^\ast_\minm=C_{d,i}$ and $\beta^\ast$ is given by \eqref{eq:CriticalRisk}.

If the energy budget is lower than $E_\minm^\ast$ 
then there does not exist a winning strategy, and we do not define $\beta_{\minm}(E^\ast)$. 
We see that for sufficiently high $E^\ast$, $\beta_{\minm}(E^\ast)=0$, so that very fit riders (who have 43\% energy than the peloton) need not assume any risks at all. 

\subsubsection{Dependence on crash propagation}
\begin{figure*}
    \centering
    \includegraphics[width=0.99\linewidth]{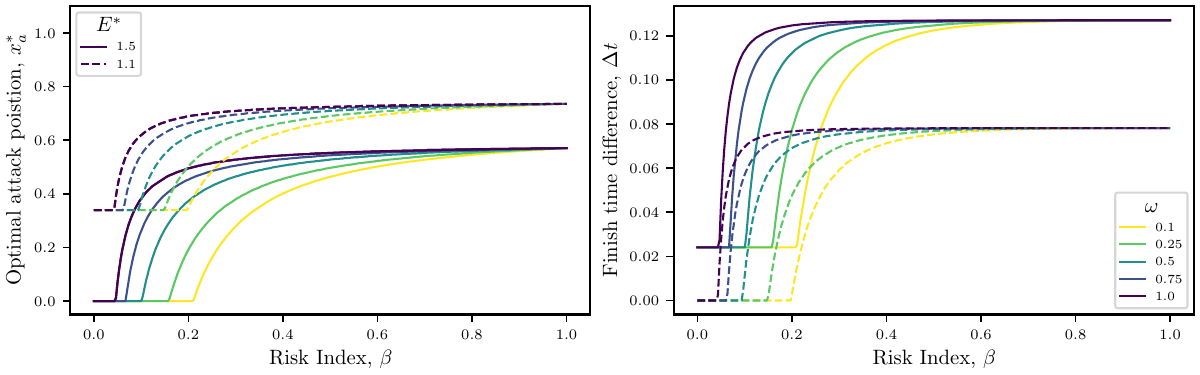}
    \caption{Optimal attack position $x_a^*(\beta)$ (left panel) and finish time ahead of peloton (right panel) plotted as a function of the risk index for $E^\ast=1.5$ (solid) and $E^\ast=1.1$ (dashed) and $\omega=0.1, 0.25,0.5,0.75, 1.0$ (denoted by colour).  The remaining parameter values are $i=5$, $C_{d,i}=0.46$, $C_{d,1}=1.43$ $\lambda=0.25$, $\mathbb{P}(C)=2$ and $N=75$.
    }
    \label{fig:crash_propagation}
\end{figure*}

In this section, we investigate how conditions such as rain, fog, road narrowing, perilous corners, or exhaustion could modify the optimal racing strategy. Mathematically, this is implemented by changing the value of the crash propagation parameter, $\omega$. Since the crash propagation is inversely proportional to $\omega$, larger values of $\omega$ may be interpreted as safer cycling conditions. In the left panel of fig.~\ref{fig:crash_propagation}, we plot the optimal attack position as a function of risk index for various crash-propagation parameters, $\omega$. When the risk and energy budget are fixed, stronger crash propagation (smaller values of $\omega$) is associated with an earlier optimal attack location (smaller $x_a^\ast$), as expected from intuition. Furthermore, we observe the minimum risk required to win the race decreases as $\omega$ increases and the riding conditions become safer.

As smaller values of $\omega$ are associated with higher crash probabilities $\mathcal{P}$ (see left panel of fig.~\ref{fig:metric_example}), stronger crash propagation manifests itself by increasing the value of the crashing component in the objective function $\mathcal{M}$, resulting in the cyclist assigning lower priority to finishing ahead of the peloton. This effect can be seen in the right panel of fig.~\ref{fig:crash_propagation}, where $\Deltat$ decreases with decreasing $\omega$, for both energy budgets studied.


\section{Towards a more realistic breakaway attempt} \label{sec:realistic_attack}

Although the model introduced in \S\ref{sec:simple_attack} can provide team directors with useful reference points for deciding team strategy, such as the minimum energy expenditure, and a good estimate for the optimal attack position, it assumes an unrealistic evolution of the power developed by a rider attempting a breakaway. Once an athlete exceeds their anaerobic limit, their body produces lactic acid, a toxin that can only be expelled in mammals by reducing the heart rate and returning to respiratory levels closer to normal~\cite{lucia_physiology_2001,joyner_endurance_2008}. 

\subsection{Modelling fatigue}

In order to incorporate the production of lactic acid, and more generally fatigue, 
we allow for a time-dependent attack power, $P_a(t)$. We follow a qualitatively similar approach to \cite{gaulOptimizingBreakawayPosition2018}, but choose instead to continue using breakaway power rather than breakaway force: 
\begin{equation}\label{eq:time_dep_power}
    P_a(t)=\begin{cases}
         P_l & t<{\ta},\\
         (P_{\maxm}-P_s)e^{-\mu(t-{\ta})}+P_s & t>{\ta}.
     \end{cases}
\end{equation}
Here, $P_{\maxm}$ is the maximum power the rider can exert, $P_s$ is the maximum sustainable power for the rider, and $P_l$ is the lurking power the rider exerts whilst in the peloton. Unless explicitly stated, we take $P_s=P_l$. The dimensionless parameter $\mu$ captures the effect of lactic acid build-up on the rider; lower values of $\mu$ correspond to a greater resistance to lactic acid. A useful interpretation for $\mu$ in the context of competitive cycling is to distinguish between riders who specialize in short, explosive outbursts of extremely high power (`sprinters') and cyclists that are better suited for endurance and long periods of sustained moderate high power, such as time trial specialists and `climbers'.

Our methodology to find the optimal breakaway position is similar to \S\ref{sec:simple_attack}, but requires additional care and numerical computation to handle the time-dependent power. We assume the rider has budgeted an energy expenditure $E^\ast$ for the race, and they wish to minimize the objective function $\mathcal{M}$ subject to some risk level $\beta$. For $P_a(t)$ as given in \eqref{eq:time_dep_power}, the total energy expenditure is
 \begin{equation}\label{eq:energy_computation_time_dep}
     E^\ast = P_l{\ta}+P_s({\tf}-{\ta})+\frac{1}{\mu}\left(1-e^{-\mu({\tf}-{\ta})}\right)(P_{\maxm}-P_s),
 \end{equation}
where ${\tf}$ is the (\emph{a priori} unknown) finish time for the rider. Given the energy budget, we may solve for the maximum attack power,
\begin{align}\label{eq:MaximumAttackPowerTimeDep}
    P_{\maxm}=P_s+\frac{\mu\re^{\mu {\tf}}\lr{E^\ast-P_s{\tf}}}{\re^{\mu {\tf}}-\re^{\mu x_a}}.
\end{align}
The velocity after the attack is 
\begin{equation}\label{eq:velocity_complicated}
    v_i(t) \equiv \frac{\rd x_i}{\rd t} =\frac{1}{C_d^{1/3}}\left(P_s +(P_{\maxm}-P_s)e^{-\mu (t-{\ta})}\right)^{1/3},\hspace{1mm} t>{\ta},
\end{equation}
using \eqref{eq:dimensionless_individual_rider} in the asymptotic limit $\varepsilon=0$. 
The position of the rider may be readily obtained by integrating \eqref{eq:velocity_complicated}:
\begin{equation}\label{eq:complicated_position}
    x_i(t)=x_a+\frac{1}{C_{d,1}^{1/3}}\int_{{\ta}}^t\left(P_s +(P_{\maxm}-P_s)e^{-\mu (t'-{\ta})}\right)^{1/3}\mathrm{d}\,t'.
\end{equation}
The finish time for the special rider, $t_f$, is obtained by solving $x_i({\tf})=1$, which must be done numerically. This provides $\Deltat=1-t_f$, and allows us to calculate the objective function $\mathcal{M}$ from \eqref{eq:metric_minima_definition}, where the crashing component is unchanged from the simpler attack \eqref{eq:CrashingProbabilitySimpleAttack}. 

Since we cannot find ${\tf}(x_a,P_s,\mu,E^\ast)$ explicitly from \eqref{eq:complicated_position}, we must minimize $\mathcal{M}$ using a numerical approach. This is achieved by solving a constrained optimization problem for three variables: the attack position $x_a$ (which is the same as the attack time ${\ta}$), the maximum attack power $P_\maxm$ and the finish time of the special rider ${\tf}$. More concretely, we minimize $\mathcal{M}(x_a, {\tf}, P_\maxm;E^\ast,\mu, P_s,\omega)$  over $x_a$, ${\tf}$, and $P_\maxm$,
subject to \eqref{eq:energy_computation_time_dep}, and \eqref{eq:complicated_position} with $x_i(t)$ set to 1. 

\subsection{Results for time-dependent attack power}
\subsubsection{Effect of fatigue}
\begin{figure*}
    \centering
    \includegraphics[width=0.99\linewidth]{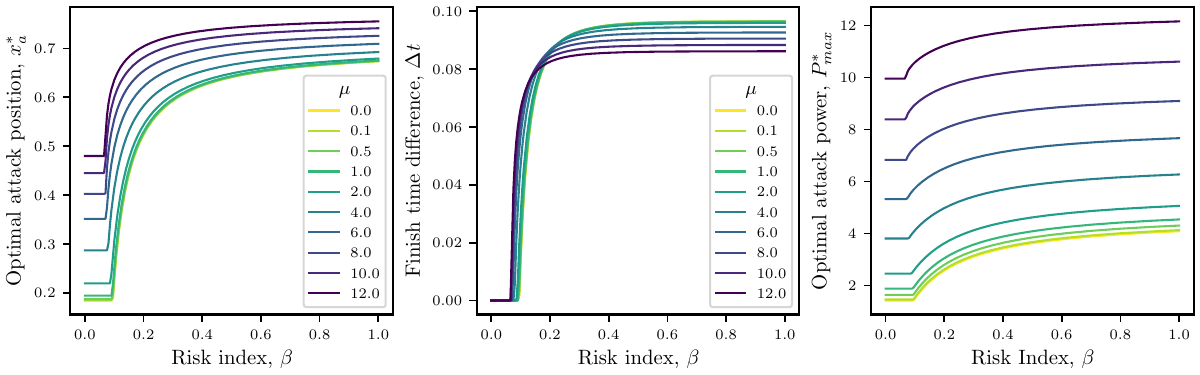}
    \caption{Left panel: optimal attack position $x^\ast_a$ plotted against risk index $\beta$, for several values of $\mu$. Middle panel: how far ahead the breakaway finishes plotted against risk index. Right panel: Optimal attack power, $P_{max}^\ast$ as a function of the risk index $\beta$. Parameter values are $E^\ast=1.25$, $i=5$, $P_s=P_l=C_{d,i}=0.46$, $C_{d,1}=1.43$, $\lambda=0.25$, $\omega=0.5$, $\mathbb{P}(C)=2$ and $N=75$.
    }
    \label{fig:mu_sweep}
\end{figure*}

In this section, we wish to understand the effect of rider fatigue on the optimal attack time, ${\ta}$, and the time difference at the finish line, $\Deltat={\tp}-{\tf}=1-{\tf}$. In the left panel of fig.~\ref{fig:mu_sweep}, we plot how the optimal breakaway time changes as a function of the risk index for different values of $\mu$. For small values of $\mu$, the breakaway power decays sufficiently slowly that it remains close to constant and we recover the solution from \S\ref{sec:simple_attack}. 

As $\mu$ increases, the optimal breakaway point moves \linebreak closer to the finish line. This is because the attack velocity eventually falls below the peloton speed (due to the increased drag associated with riding without drafting). Thus, attacks that take place early in the race run the risk of failing in the final stages of the race due to the exhausted rider being caught by the peloton. Hence, we expect the optimal attack location to increase as $\mu$ increases (see left panel of fig.~\ref{fig:mu_sweep}).

In the middle panel of fig.~\ref{fig:mu_sweep} we plot how the time difference at the finish line changes with risk index if the breakaway is attempted at the optimal attack position $x^\ast_a$ (from the left panel). 
For moderate to high risk, we observe that endurance riders (small $\mu$) finish the race with a larger margin than explosive riders (large $\mu$). Again, as $\mu\rightarrow 0$ we recover the solution from the constant-breakaway-power model, plotted in yellow. 
As in the simpler model from \S\ref{sec:simple_attack}, there is a minimum risk that must be assumed to finish ahead of the peloton. Notably, this required minimum risk decreases as $\mu$ increases. However, consistent with expectations, the time gained over the peloton also diminishes with increasing $\mu$. Practically speaking, this implies that cyclists with greater endurance capacity are better positioned to establish and maintain substantial leads over the peloton. On the other hand, explosive cyclists attack closer to the finish line and reach the finish line with the peloton close behind.

In the right panel of fig.~\ref{fig:mu_sweep} we plot how the optimal maximum attack power changes with the risk index. As expected, larger values of $\mu$ lead to larger maximum attack powers, and for $\mu\gg1$ we we see that $P_\maxm=\mathcal{O}(\mu)$, a scaling predicted by \eqref{eq:MaximumAttackPowerTimeDep}. As $\mu \rightarrow 0$, we recover the $\mu=0$ trends from fig.~\ref{fig:x_optimal_vs_risk}. We remark that the large values of $\mu$ we consider here ($\mu\approx 10$) require the generation of high attack powers $P_\maxm\approx10$. In dimensional terms, this corresponds to exerting a peak power output that is 10 times higher than the peloton average during the start of the attack. This is not too dissimilar to real data for elite sprinters, who have been reported to attain peak powers in excess of 1000 W during bursts \cite{kordi_mechanical_2020}. This estimate thus compares well with peloton powers for flat stages of around 150 W~\cite{bernalStrava2021}.

\subsubsection{Role of energy budget}
\begin{figure*}
    \centering
    \includegraphics[width=0.99\linewidth]{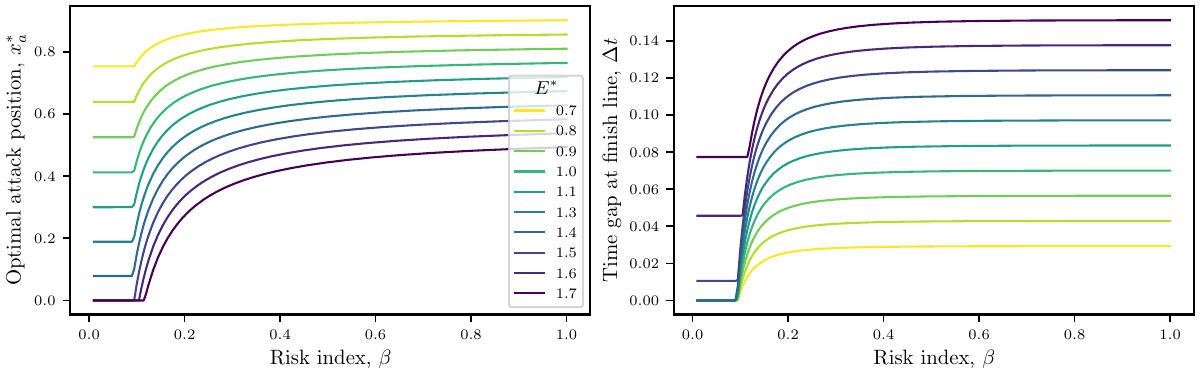}
    \caption{Left panel: optimal attack position $x^\ast_a$ plotted against risk index $\beta$, for several values of the energy budget $E^\ast$. Right panel: how far ahead the breakaway finishes plotted against risk index. 
    Parameter values are $\mu=1$, $i=5$, $P_s=P_l=C_{d,i}=0.46$, $C_{d,1}=1.43$, $\omega=0.5$, $\lambda=0.25$, $\mathbb{P}(C)=2$ and $N=75$.}
    \label{fig:E_sweep}
\end{figure*}

We look at how the energy budget changes the optimal breakaway strategy in fig.~\ref{fig:E_sweep}. We fix the fatigue parameter to $\mu=1$, an intermediate value where the rider gets exhausted in a time-scale comparable to the time to finish the race. In the left panel, we plot how the optimal attack position changes with the risk index $\beta$, for several values of the energy budget $E^\ast$. As in the simpler model from \S\ref{sec:simple_attack}, for large values of the energy budget, the optimal solution for suitably risk-averse cyclists (small $\beta$) is to break away at the start of the race. For risk indices higher than the minimum risk, the optimal attack position gradually increases, approaching a an energy dependent limit as $\beta\rightarrow1$. In the right panel of fig.~\ref{fig:E_sweep}, we plot how far ahead the peloton the breakaway finishes as function of the risk index for several values of $E^\ast$. As expected from the previous model, small values of $E^\ast$ require the cyclist to assume a minimum risk to win.


\section{Discussion and future directions} \label{sec:discussion}

In this paper, we have developed a mathematical framework that captures the key trade-offs faced by elite cyclists: balancing aerodynamic efficiency, energy expenditure, and the risk of crashing. By incorporating both physiological limits and probabilistic crash modelling, our approach extends earlier models from \cite{gaulOptimizingBreakawayPosition2018} to more realistically capture the strategic considerations involved in breakaway attempts.

One of the key strengths of our model lies in its relative simplicity and analytical tractability. Despite this, it successfully quantifies several intuitive yet previously unmodelled strategic insights. For example, it confirms that breakaways near the end of a race are more energy-efficient, while earlier attacks can result in larger time gains, but at a higher energetic cost. Similarly, spending more time near the rear of the peloton reduces drag~\cite{blockenAerodynamicDragCycling2018} but significantly increases the risk of being caught in a crash. These results not only validate widely held beliefs among athletes and team directors but also allow for rigorous quantitative evaluation of such strategies.

A particularly novel feature of this work is the inclusion of crash dynamics as a fundamental part of strategic planning. By treating the probability of crashing as a dynamic penalty, we show that risk management is as crucial to success as fitness or power output. This aligns with real-world race scenarios, where riders must weigh the dangers of remaining in a densely packed peloton against the energetic demands of escaping it. In elite competition, where physical capabilities are often comparable among top athletes, it is these marginal strategic decisions -- timing of attacks, positioning, and risk tolerance -- that frequently determine the outcome.

The model is also versatile enough to accommodate different risk profiles via the parameter $\beta$, representing a rider’s willingness to trade risk for performance. This provides a useful lens to interpret decisions made by different types of riders: for instance, a General Classification (GC) leader may choose a conservative strategy, while a \emph{stage hunter} might favour risky strategies.

While the current model assumes a flat course and relatively simple power profiles, extensions are straightforward. In Appendix~\ref{sec:mountains}, we outline a framework for incorporating elevation changes, showing that our methods readily generalize to more complex terrains. Though analytical results are harder to obtain in this setting, standard numerical solvers can be used to explore optimal strategies in varied topography. In Appendix~\ref{sec:complex_crash}, we show how the crash model could be generalized to account for more realistic crash propagation patterns and environment-dependent risks such as road conditions, weather, or race congestion.

There are several promising directions for future work. One is to include dynamic chasing behaviour, where the peloton adjusts its speed in response to breakaways, 
and how this affects the viability of a breakaway. Another is to model physiological decay more precisely by linking the power decay rate $\mu$ to underlying rider attributes or empirical data (e.g., from power meters or lactate threshold tests~\cite{jobson_analysis_2009}). Similarly, incorporating nutrition and hydration strategies, particularly over multi-stage races such as the Tour de France, could provide insights into optimal energy allocation over days rather than single stages \cite{earnest_relation_2004}.

Another future direction for this work could be to explore velocity-dependent crash probabilities, which would reflect reduced reaction times at higher speeds. Moreover, integrating spatial features of the course such as  narrow roads, roundabouts, or downhill corners would offer a more nuanced crash model. These refinements could also be informed by data analytics from race telemetry and crash reports.

A third future direction is to understand how head-, tail- and cross-winds can shape racing strategies, particularly when the weather forecast is different for different stages of a race. This is because strong headwinds should lead riders to prioritize staying inside the peloton due to higher drag at the same ground speed, but tail winds effectively reduce drag and hence the advantages of staying in the peloton. Moreover, crosswinds can lead to the formation of echelons \cite{beaumont_fighting_2023}, where riders who are not well positioned at the front of the peloton may be easily dropped \cite{phillips_determinants_2020}, giving rise to exhilarating racing conditions.

Finally, we emphasize that elite racing is not solely a test of physical ability, but of strategic timing and decision-making. As athletic performance nears physiological limits~\cite{joyner_endurance_2008,haake_impact_2009}, small tactical choices, such as when and where to initiate a breakaway, can determine the outcome of a race. Our framework provides a structured way to analyse these choices, translating race-day intuition into quantifiable insights. The analytical clarity that our  mathematical model provides makes this a valuable tool for interpreting rider behaviour and guiding tactical planning in competitive settings.

\textbf{Acknowledgments:} We are grateful to J.~W.~M.~Bush for his feedback on one of the initial ideas for this project and to Alberto Conde Mellado and Lucas T.~Y.~Yeung for helpful discussion. The photograph in this text belongs to Luis \'{A}ngel G\'{o}mez, to whom we thank for his generosity. For the purpose of Open
Access, the authors will apply a CC BY public copyright license to any Author
Accepted Manuscript version arising from this submission.

\textbf{Funding:} JCV is funded by a St. John's College, Oxford scholarship. 



\appendix

\section{An in-depth look at the attack microstructure}\label{sec:transition_layers}

\begin{figure}
    \centering
    \includegraphics[width=0.99\linewidth]{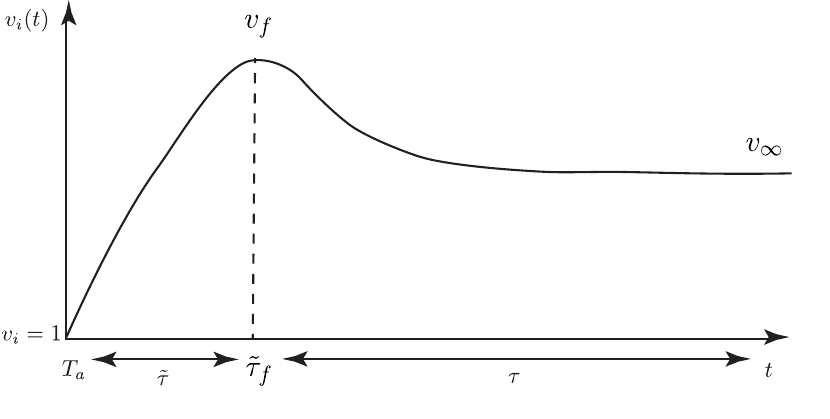}
    \caption{A sketch of the asymptotic transition layers during the early stages of an attack attempt.
    }
    \label{fig:transition_layers_sketch}
\end{figure}

In our analysis in \S\ref{sec:scalings}, we exploit the fact that $\epsilon \ll 1$ to obtain a leading-order solution for the trajectory of the rider by setting $\epsilon=0$. In doing so, we solve \eqref{eq:dimensionless_peloton} and \eqref{eq:dimensionless_individual_rider} in such a way that the sudden jump in power instantaneously translates to a jump in speed. In reality, there is a small time when the rider is accelerating from the peloton speed to their attack speed, as fig.~\ref{fig:transition_layers_sketch} shows. There, we numerically solve the nonlinear problem \eqref{eq:dimensionless_individual_rider} for small but finite $\epsilon =0.005$. The velocity first rapidly increases when the attack occurs, as the rider benefits from accelerating with reduced drag due to the presence of the peloton. Once the rider reaches the head of the peloton (and the breakaway truly commences) the rider experiences increased drag, with their speed reducing until stabilizing to the attack velocity computed in the main text, $v_a=(P_a/C_{d,1})^{1/3}$. We can understand this transition region by examining the time close to the attack position.

For a typical race length of $150$\,km, $\delta \approx 2\times 10^{-5}=\gamma\epsilon^2$, where $\gamma \approx 1$. This means that the approach to the equilibrium velocity after the attack takes place on two timescales: the first, O$(\epsilon^2)$,  timescale captures the motion of the rider to the front of the peloton; the second, $O(\epsilon)$,  timescale captures the part of the race once the rider has escaped the peloton and their speed approaches the solo equilibrium value. 

To capture the first phase, we scale $t={\ta}+\epsilon^2 \tilde{\tau}$ and $x_i=x_a+\delta\zeta_i$, where $\delta=\gamma\epsilon^2$. Substituting these scalings into \eqref{eq:dimensionless_individual_rider} gives, to leading order in $\epsilon$: 
\begin{align}
\gamma m_i\frac{\rd^2\zeta_i}{\rd \tilde{\tau}^2}=P_i-C_{d,i}(\zeta_i),
\end{align}
where 
\begin{align}
    C_{d,i}(\zeta_i)=\frac{\hat{C}_d^{\minm}+\lr{\hat{C}_d^{\maxm}-\hat{C}_d^{\minm}}\re^{-\lambda\zeta_i}}{\langle \hat{C}_d\rangle}.
\end{align}
We solve this subject to $\zeta_i=-(i-1)$ and $\rd \zeta_i/\rd \tilde{\tau}=0$ for $\tau=0$, which provides the solution until $\tilde{\tau}=\tilde{\tau}_d$ defined by $\zeta_i(\tilde{\tau}_d)=0$, which corresponds to the time when the rider has reached the front of the peloton. Note that here we are solving for $\zeta_i$ as a continuous variable. This allows us to track the motion of the special rider past the discrete values $\zeta_i=i-1$ that correspond to locations at which we permit riders to sit while riding in the peloton and not making an attack. The speed of the rider when they reach the front of the peloton is $v_f \equiv\rd x_i/\rd t\vert_{\tilde{\tau}=\tilde{\tau}_f}=1+\rd \zeta_i(\tau_f)/\rd \tilde{\tau}.$ 

For the second part of this, the appropriate distinguished limit is obtained by scaling $t={\ta}+\epsilon\tau$ and $x_i=x_a+\epsilon X_i$. Substituting these scaling into \eqref{eq:dimensionless_individual_rider} gives, to leading order in $\epsilon$: 
\begin{align}
\label{eq:microstructure second timescale}
    m_i\frac{\rd^2 X_i}{\rd \tau^2}=P_i\lr{1+\frac{\rd X_i}{\rd\tau}}^{-1}-C_{d,i}(0)\lr{1+\frac{\rd X_i}{\rd \tau}}^2.
\end{align}
We solve this subject to $X_i=0$ and $\rd X_i/\rd \tau=v_f-1$ at $\tau=0$, which ensure that the solution matches to the $O(\epsilon^2)$ timescale of the rider moving to the front of the peloton. We note that we can write \eqref{eq:microstructure second timescale} in terms of the velocity $\rd X_i/\rd \tau=\rd x_i/\rd t-1=v_i$, which may then be solved straightforwardly subject to $v_i=v_f$ when $\tau=0$, to give   
\begin{align}
v_i=\frac{P_i v_f}{C_d(0)v_f+(P_i-C_d(0)v_f)\re^{-P_i \tau/m_i}}.
\end{align}
The composite solution for $v_i$ first rises as the rider makes their way through the peloton, reaching a maximum when they reach the front of the peloton, before decaying exponentially to $v_\infty \equiv P_i/C_d(0)$ once the rider escapes the peloton and becomes a solo rider (see fig.~\ref{fig:transition_layers_sketch}).  
\color{black}

\section{Minimum risk computation}\label{sec:minimum_risk_appendix}


To compute the value of the minimum risk, 
 we first recall the minimum viable attack position $x_a^\minm$ is given by \eqref{eq:minimum_attack_position}.
For $E^\ast_\minm<E^\ast<E^\ast_{\mathrm{crit}}$, the minimum attack position lies in the open interval $x_a^\minm\in(0,1)$. We focus on this range of energy budgets as (i) for $E^\ast>E_\mathrm{crit}^\ast$ the start location $x_a=0$ is a viable attack position and conservative riders will prefer it, and (ii) for $E^\ast<E_\minm^\ast$ there is no viable attack position (mathematically, $x_a\geq1$).
Moreover, for $E^\ast_\minm<E^\ast<E^\ast_{\mathrm{crit}}$ the objective function has a discontinuity in its derivative at $x_a^\minm$ for all values of $\beta$. In this case, as we can see in the right panel of fig.~\ref{fig:metric_example}, for $x_a<x_a^\minm$ the objective function is flat. For $x_a>x_a^\minm$, the behaviour of the objective function depends on the risk index, and there are two qualitatively distinct regimes. 

The first regime occurs for small $\beta$, and is characterized by $\mathcal{M}(x_a)$ being a monotonically increasing function for $x_a>x_a^\minm$. In this case, $x_a^\dagger$ is not defined as there is no point where $\mathcal{M}'(x_a)=0$ in the domain (see, for example, the curves for $\beta=0$ and 0.1 in the right panel of fig.~\ref{fig:metric_example}). In this regime, $x_a\ast=x_a^\minm$.

A second regime exists for larger $\beta$, where $\mathcal{M}(x_a)$ is a non-monotonic function of $x_a>x_a^\minm$ (see, for example, the curves for $\beta=0.4$, 0.6, 0.8 and 1.0 in the right panel of fig.~\ref{fig:metric_example}). In this case, $x_a^\dagger$ is well defined as the point where $\mathcal{M}'(x_a^\dagger)=0$, and this is the global minimum of $\mathcal{M}$, so that $x_a^\ast=x_a^\dagger$.

The transition between both regimes occurs at a critical $\beta$, which we define as $\beta^\ast$. To find this value, we observe that in the second regime we have $\mathcal{M}'(x_a^\minm)\leq0$, while in the first regime we have $\mathcal{M}'(x_a^\minm)\geq0$. At the critical value of $\beta^\ast$, we are in an intermediate regime that satisfies exactly $\mathcal{M}'(x_a^\minm)=0$. An alternative way of seeing this is to consider how $x_a^\dagger$ changes as $\beta$ is decreased from a high value. As $\beta$ is decreased, $x_a^\dagger$ decreases, eventually reaching $x_a^\dagger=x_a^\minm$ for $\beta=\beta^\ast$. By definition $x_a^\ast$ is the point where the derivative vanishes, so that substituting we also find $\mathcal{M}'(x_a^\minm)=0$. Using the definition of $\mathcal{M}$, \eqref{eq:metric_minima_definition}, this condition reads
\begin{equation}
    -\beta^\ast\left.\frac{\partial\Deltat}{\partial x_a}\right\vert_{x_a=x_a^\minm}+(1-\beta^\ast)\left.\frac{\partial \mathcal{P}}{\partial x_a}\right\vert_{x_a=x_a^\minm}=0,
\end{equation}
and solving for the critical risk $\beta^\ast$,
\begin{equation}
    \beta^\ast=\frac{\left.\frac{\partial \mathcal{P}}{\partial x_a}\right\vert_{x_a=x_a^\minm}}{\left.\frac{\partial \mathcal{P}}{\partial x_a}\right\vert_{x_a=x_a^\minm}+\left.\frac{\partial\Deltat}{\partial x_a}\right\vert_{x_a=x_a^\minm}}.
\end{equation}
We can compute the partial derivatives of the objective function components. For the time difference, 
\begin{equation}
    \frac{\partial \Deltat}{\partial x_a}=1+\frac{3}{2}\left(\frac{1-x_a}{E^\ast-C_{d,i}}\right)^{1/2}C_{d,1}^{1/2}-\frac{C_{d,i}C_{d,1}^{1/2}}{2}\left(\frac{1-x_a}{E^\ast-C_{d,i}}\right)^{3/2}.
\end{equation}
Now, evaluating $\partial \Deltat/\partial x_a$ at $x_a=x_a^\minm$, the above simplifies to
\begin{equation}
\left.\frac{\partial \mathcal{P}}{\partial x_a}\right\vert_{x_a=x_a^\minm}=\frac{1}{2}\left(1-\frac{C_{d,i}}{C_{d,1}}\right),
\end{equation}
where we have used the fact that 
\begin{equation}
    \eta_\minm^2 \equiv \frac{1-x_a^\minm}{E^\ast-C_{d,i}x_a^\minm}=\frac{1}{C_{d,1}}.
\end{equation}
Since the probability of crashing is a linear function of $x_a$ from \eqref{eq:CrashingProbabilitySimpleAttack}, its derivative is constant and given by
\begin{equation}
    \left.\frac{\partial\mathcal{P}}{\partial x_a}\right\vert_{x_a=x_a^\minm}=\frac{\partial\mathcal{P}}{\partial x_a}=\frac{\mathbb{P}(C)}{N}\left(\frac{1-e^{-\omega i}}{1-e^{-\omega}}-1\right).
\end{equation}
Thus, the critical risk is
\begin{equation}
    \beta^\ast=\frac{\frac{\mathbb{P}(C)}{N}\left(\frac{1-e^{-\omega i}}{1-e^{-\omega}}-1\right)}{\frac{\mathbb{P}(C)}{N}\left(\frac{1-e^{-\omega i}}{1-e^{-\omega}}-1\right)+\frac{1}{2}\left(1-\frac{C_{d,i}}{C_{d,1}}\right)}.
\end{equation}

\section{Elevation}\label{sec:mountains}

In this section, we outline how to generalize the methods from the main paper to account for elevation changes, including an example simulation of a race over arbitrary terrain. In this case, the generalized equation for the peloton motion, \eqref{eq:dimensionless_peloton}, when riding on a course of arbitrary steepness $\theta(x)$, is given by
\begin{equation}\label{eq:dimensionless_peloton2}
    \epsilon \frac{\rd^2 x_p }{\rd t^2}=P(t)\left( \frac{\mathrm{d}x_p}{\mathrm{d}t}\right)^{-1}-\left(\frac{\mathrm{d}x_p}{\mathrm{d}t}\right)^2-\gamma\sin(\theta(x_p))=0,
\end{equation}
where 
\begin{align}
    \gamma &= \frac{2^{1/3}\langle \hat{m}\rangle g }{\hat{P}_p^{2/3}(\langle \hat{C}_d\rangle \rho A)^{1/3}}\approx 40,
\end{align}
which is a ratio between the force due to gravity and aerodynamic drag. Although $\gamma$ appears rather large, since the inclination angle $\theta$ is at most around 0.1, the effect of gravity is $\mathcal{O}(1)$. The generalized equation for motion of an individual rider, \eqref{eq:dimensionless_peloton}, is 
\begin{equation}
    \epsilon m_i \frac{\rd^2 x_i}{\rd t^2}={P_i(t)}\left(\frac{\rd{x}_i}{\rd t}\right)^{-1}-C_{d,i}\left(\frac{\rd x_i}{\rd t}\right)^2-m_i\gamma\sin(\theta(x_i))=0,\label{eq:dimensionless_individual_rider with elevation} 
\end{equation}
One crucial difference with the flat course presented in the main paper is that the time it takes for the peloton to complete the course, ${\tp}$, is no longer unity, as it will depend on the race profile $\theta(x)$. As a result, the total energy expenditure in the peloton is no longer unity either, but is now given by $E={\tp}$, using \eqref{eq:dimensionless energy}. 

We now simulate a simple attack, where the rider is located at some position in the interior of the peloton until they reach $x_a$, where they attempt to breakaway with a constant power, as in \S\ref{sec:simple_attack}. Unlike for a flat course, analytical progress is complicated, so we solve for the trajectory of the peloton and the rider numerically. The solution strategy for $x_p(t)$ and $x_i(t)$ is: 
\begin{itemize}
    \item We evolve the position of the peloton with a constant power $P_p(t)=1$, solving \eqref{eq:dimensionless_peloton2}. 
    \item Before the attack, for $x_i(t)<x_a$ the rider simply moves at the same speed as the peloton, $\dot{x}_i=\dot{x}_p$, obtained in the first step. Using \eqref{eq:dimensionless_individual_rider with elevation} we can find the power $P_i(t)$ that the rider must exert to maintain this velocity, which will be in general lower than 1, as the rider is located inside of the peloton and experiences lower drag. 
    \item At $x_i(t)=x_a$, the rider attacks so that we now solve \eqref{eq:dimensionless_individual_rider with elevation} with the power set to $P_i(t)=P_a$ (this could also be an exponentially decaying power profile, as in \S\ref{sec:realistic_attack}). 
    \item From here on we solve the (now uncoupled) equations of motion for $x_p(t)$ and $x_i(t)$, until both of them reach the end of the race at $x=1$, where an event handler automatically stops the integration.
\end{itemize}

We use Scipy's solve\_ivp routine to solve the problem numerically. The routine can handle ``events", stopping the integration once a certain condition has been met, namely reaching the finish line at $x=1$. The integration method is either RK45, or BDF if we suspect the sharp power profile leads to a stiff problem. 

\begin{figure*}
    \centering
    \includegraphics[width=0.99\linewidth]{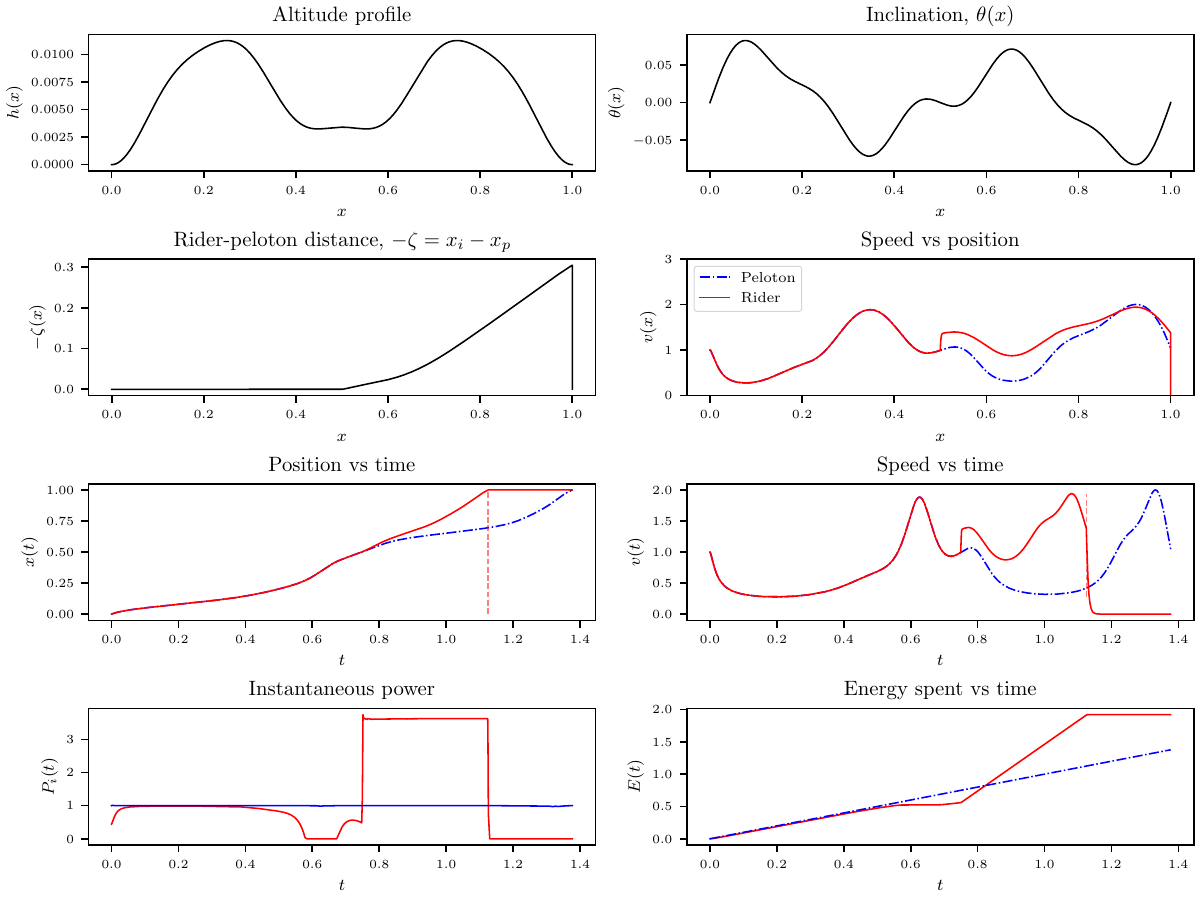}
    \caption{Details for the simulation of a rider (solid red lines) breaking away from the peloton (dashed blue lines) at $x_a=0.5$, in a hilly race.The parameters are $i=5$, $C_{d,1}=1.43$, $P_a=3.6$, $x_a=0.5$, $N=75$ and $\lambda=0.25$. 
    The feint red line in row 3, left  indicates the point where the individual rider has finished the course. 
    }
    \label{fig:elevation_example}
\end{figure*}

In fig.~\ref{fig:elevation_example} we show the results for a simulation of a single rider attacking at position $x_a=0.5$, with a constant attack power as in \S\ref{sec:simple_attack}. In the top row we can see the example mountainous course elevation profile, $h(x)$, obtained by a superposition of sines and cosines. The elevation is small because it is scaled by the course length $L$. The steepness is given by $\theta(x)=\tan^{-1}(h'(x))$. In the second row, we plot how the distance between the rider and the peloton widens from $x=x_a$ onwards, as well as the velocity of the rider (red) and the peloton (blue) as a function of position. We see that the individual rider moves faster in flat (and uphill) sectors of the course, while the peloton descends faster on the steepest regions due to the reduced drag it experiences. The third row shows the positions and velocities of the rider and the peloton, but now as a function of time. The bottom row shows  the evolution of the power of the rider, $P_i(t)$, which varies in such a way that the velocity $v_i=1$ before the breakaway, and $v_i=P_a$ after the breakaway. We do not allow for negative powers, so that in very steep downhill sections where the rider needs to do no work to keep up with the peloton, the power is set to zero. The peloton power is set to unity for the whole course. The lines in the right panel of the bottom row correspond to the cumulative energy consumed, defined by \eqref{eq:dimensionless energy}.

\section{Complex crash dynamics}\label{sec:complex_crash}
\subfile{sections/appendix_complex_crash.tex}

\bibliographystyle{ieeetr}
\bibliography{cycling_project}

\end{document}

%% file: sections/appendix_complex_crash.tex
If we drop the assumption that the distribution of the position of crash start is uniform, more complicated crashing dynamics emerge. We recall that the quantity of interest is the probability of being involved in a crash conditioned on a crash occurring, which we define as $\mathbb{P}(C_j\vert C)$. Let $\mathbb{P}(S_k\vert C)=s_k$ be the probability that the cyclist in position $k$ crashes first, which we now allow to vary with $k$. Once a crash has been initiated at position $k$, the probability that the cyclist in position $j$ crashes too is zero if $j<k$ (rider is ahead of the crash), and will in general depend on the distance between $j$ and $k$ for $j>k$. Therefore, we write 
\begin{equation}
    \mathbb{P}(C_j\vert S_k)=\begin{cases}
        0 & j<k \\
        f(j-k) & j\geq k,
    \end{cases}
\end{equation}
where $f(\cdot)$ is a decreasing function of the distance between $j$ and $k$ (the further away the less likely to also crash) and satisfies $f(0)=1$ ($j=k$ is the cyclist that starts the crash, so they crash with probability 1). In general, both the distribution of the crash start position $s_k$ and the crash propagation statistics could be learned from real-life data. In the main paper, we used a uniform distribution for $s_k=1/N$ and an exponential decay for $f(j-k)=e^{-\omega(j-k)}$ as realistic yet simple starting points. Furthermore, the parameters and distributions will in general depend on the speed of the peloton (which will affect reaction time), environmental factors (for example, fog, or wet roads that decrease braking performance) and how densely packed the peloton is.

Once their function or numerical form has been found, the probability that the cyclist at position $j$ crashes given that a crash has occurred can be found using the law of total probability:
\begin{equation}
    \mathbb{P}(C_j\vert C)=\sum_{k=0}^N\mathbb{P}(C_j\vert S_k)\mathbb{P}(S_k)=\sum_{k=0}^jf(j-k)s_k.
\end{equation}